\documentclass[11pt,reqno]{amsart}
\usepackage{amsmath}
\usepackage{amssymb}
\usepackage{amsthm}
\usepackage{mathtools}
\usepackage{color}
\usepackage{eucal}
\usepackage{tikz}
\usepackage{gastex}
\usepackage{stmaryrd}
\usepackage{caption,subcaption}
\usepackage{latexsym}
\usepackage{indentfirst}
\usepackage{graphicx,accents}
\usepackage{mathrsfs}

\usetikzlibrary{snakes}

\DeclareSymbolFont{rsfscript}{OMS}{rsfs}{m}{n}
\DeclareSymbolFontAlphabet{\mathrsfs}{rsfscript}

\numberwithin{equation}{section}
\newtheorem{prop}{Proposition}[section]
\newtheorem{teor}[prop]{Theorem}
\newtheorem{lem}[prop]{Lemma}
\newtheorem{cor}[prop]{Corollary}

\def\Jc{\mathrel{\mathrsfs{J}}}
\def\Dc{\mathrel{\mathrsfs{D}}}
\def\Hc{\mathrel{\mathrsfs{H}}}
\def\Lc{\mathrel{\mathrsfs{L}}}
\def\Rc{\mathrel{\mathrsfs{R}}}
\def\Kc{\mathrel{\mathrsfs{K}}}
\def\Rcc{\mathrm{R}}

\def\Lcc{\mathrm{L}}
\def\Dcc{\mathrm{D}}

\def\Kcc{\mathrm{K}}

\def\ep{\epsilon} 
\def\up{\upsilon} 
\def\la{\lambda} 
\def\ka{\kappa} 
\def\be{\beta} 
\def\si{\sigma} 

\def\ap{\approx}
\def\xr{\xrightarrow}
\def\cev{\overset{{}_{\shortleftarrow}}}

\def\ol{\overline}
\def\l{\mathrm{l}}

\newcommand{\oX}{\ol{X}}

\renewcommand{\iff}{if and only if}

\DeclareMathOperator{\ft}{FI}
\DeclareMathOperator{\fw}{F}

\DeclareMathOperator{\fo}{F_1}

\DeclareMathOperator{\La}{L}
\DeclareMathOperator{\mrx}{MR}
\DeclareMathOperator{\mx}{M}

\DeclareMathOperator{\G}{G}

\DeclareMathOperator{\R}{R}

\def\ftx{\ft(X)}
\def\fx{\fw(X)}
\def\fox{\fw^1(X)}

\def\Gd{\G_d}
\def\Gf{\G^5}
\def\Gid{\G_{i,d}}
\def\Ge{\G_e}
\def\Gie{\G_{i,e}}
\def\Gi{\G_i}

\def\Gp{\G^+}

\def\mxm{\mx^1}
\def\lx{\La}
\def\lo{\La_1}
\def\lop{\La_1^+}
\def\lom{\La_1^-}
\def\ltp{\La_2^+}
\def\ltm{\La_2^-}
\def\lt{\La_2}

\def\pres{\langle\G,\R\rangle}

\title{Weakly free regular semigroups I.\\ the general case}
\author{Lu\'\i s Oliveira}
\address{CMUP, Departamento de Matem\'atica,
Faculdade de Ci\^encias, Universidade do Porto,
R. Campo Alegre s/n, 4169-007 Porto, Portugal}
\email{loliveir@fc.up.pt}

\begin{document}

\begin{abstract}
We prove the existence of a regular semigroup $\fx$ weakly generated by $X$ such that all other regular semigroups weakly generated by $X$ are homomorphic images of $\fx$. The semigroup $\fx$ is introduced by a presentation and the word problem for that presentation is solved.
\end{abstract} 


\maketitle

\section{Introduction}

Let $S$ be a semigroup. As usual, $E(S)$ denotes the set of idempotents of $S$ and $V(x)$ denotes the set of inverses of $x\in S$, that is, the set of all elements $x'\in S$ such that $xx'x=x$ and $x'xx'=x'$. Then $S$ is regular if all sets $V(x)$ are nonempty. The property of being regular is easily lost for taking subsemigroups. In other words, if $S$ is regular and $X\subseteq S$, the subsemigroup generated by $X$ is not be regular very often. 

A regular subsemigroup $T$ of $S$ is \emph{weakly generated} by a subset $X$ of $S$ if $T$ has no proper regular subsemigroup containing $X$. In particular, $S$ is weakly generated by $X$ if $S$ has no proper regular subsemigroup containing $X$. We reinforce that a regular [sub]semigroup weakly generated by $X$ needs not to be generated by $X$ since $X$ may generate a proper non-regular subsemigroup. We also point out that the same set $X$ may weakly generate several distinct regular subsemigroups of $S$.

Since subsemigroups of regular semigroups may not be regular, we use the concept of e-variety for dealing with classes of regular semigroups instead of the more common notion of variety of algebras. An e-variety of regular semigroups is a class of these semigroups closed for regular subsemigroups, homomorphic images and direct products \cite{ha1,ks1}. This concept was intensively studied in the nineties. The key ingredient in those studies was the notion of \emph{bifree object} introduced in \cite{ks1} (we refer the reader to that paper for a concrete definition of `bifree object'). This concept was the natural adaptation for e-varieties of the more common notion of free objects in varieties of algebras. Most paper at that time studied the structure of the bifree objects in e-varieties that have them, and tried to obtain Birkhoff type theorem for e-varieties \cite{au94,au95,au95b,as,br,ks1,ks2,yeh} (see also \cite{bi,LO04,LO17,LO21}).

Unfortunately, not all e-varieties have bifree objects. Some attempts have been made to generalize the notion of bifree objects to other e-varieties with the aim of obtaining Birkhoff type theorems for those e-varieties \cite{ct,ka}. However, those attempts did not completely fulfill their goal, the goal of obtaining a theory that would allow us to work with e-varieties in a manner similarly to how we work with varieties of algebras. The fact of not existing a good notion of `free object' for all e-varieties seems to have refrain the interest in these classes. We refer the reader to \cite{au02, jo, tr} for an overview about the theory developed around the notion of e-variety of regular semigroups. See also \cite{sz20} for an overview about the general theory of regular semigroups. We should also mention that, recently, Higgins and Jackson \cite{hj} developed a theory based on `equation systems' to described elementary classes of algebras (with a specific signature) closed for homomorphic images and direct products. E-varieties of regular semigroups can be seen as special cases of these classes.

Yeh \cite{yeh} proved that only the e-varieties of locally inverse semigroups or the e-varieties of regular $E$-solid semigroups have bifree objects on any set $X$. A close analysis to the proof of Yeh's result reveals that bifree objects on any set $X$ exist in an e-variety if and only if any semigroup $S$ of that e-variety has a unique regular semigroup weakly generated by $X$, for any `matched' subset $X$ of $S$. Recall that a matched subset of a regular semigroup is a subset $A$ where all elements of $A$ have inverses inside $A$. Thus, there is a close relation between bifree objects and weakly generated regular subsemigroups. In this paper, we prove that, although the e-variety of all regular semigroups does not have bifree objects, it has objects with some universal properties related to the concept of weakly generated regular semigroups.
 
Recently, we began studying weakly generated regular semigroups. In \cite{LO22} we proved the existence (and uniqueness up to isomorphism) of a regular semigroup $\ftx$ weakly generated by $X$ idempotents such that all other regular semigroups weakly generated by $X$ are homomorphic images of $\ftx$. In \cite{LO23} we proved the existence (and uniqueness up to isomorphism) of a regular semigroup $\fo$ weakly generated by one (non-idempotent) element $x$ such that all other regular semigroups weakly generated by $x$ are homomorphic images of $\fo$. In the present paper, we generalize the main result of \cite{LO23} to any set $X$ and extend the main result of \cite{LO22} to the non-idempotent case. In other words, the main result of this paper is the proof of the existence (and uniqueness up to isomorphism) of a semigroup $\fx$ weakly generated by a set $X$ such that all other regular semigroups weakly generated by $X$ are homomorphic images of $\fx$.

This paper will follow closely the structure of \cite{LO23}. The results of \cite{LO23} are stated for $|X|=1$, but the size of $X$ is in fact irrelevant for many of their proofs. Thus, we will refrain ourselves to include here again those proofs in order to keep the paper short. We will, however, make a reference to the corresponding result of \cite{LO23}. Thus, if there is a reference to a result of \cite{LO23} in the statement of a result here, that means that the size of $X$ has no influence in the proof presented in \cite{LO23}, and so that proof is valid also for the result stated here. 

This paper is organized as follows. In the next section, we recall some basic concepts needed for this paper. In Section 3 we construct the semigroup $\fx$. This semigroups will be defined by a presentation. Since both the set of generators and the set of relations of this presentation have a complex definition, we divide this section in subsections. In Section 4 we solve the word problem for this presentation. Once again we divide this section into subsections that indicate the main steps towards the solution of the word problem. We believe that the division of these two sections in subsections will help the reader follow the paper. Section 5 is devoted to the structure of $\fx$. In this section we prove that $\fx$ is a regular semigroup weakly generated by $X$ and we describe, for example, its Green's relations (but not only). The study of the structure of $\fx$ is important to prove, in Section 6, the main result of this paper already mentioned above. In this Section 6 we prove also that all semigroups with the same universal property as $\fx$ are isomorphic to $\fx$.

\section{Preliminaries}\label{sec2}

Let $X$ be a nonempty set. The free semigroup on $X$ is the set of all (finite non-empty) sequences of elements of $X$ endowed with the concatenation operation. As usual, we denote by $X^+$ this free semigroup, and call its elements \emph{words} (on $X$). The elements of $X$ are also designated by \emph{letters}. 

The \emph{length} $\l(u)$ of a word $u=x_1x_2\cdots x_n$, with $x_i\in X$ for $1\leq i\leq n$, is its number of letters, that is, $n$ for this case; while the \emph{content} of $u$ is the set of distinct letters from $X$ that occur in $u$. We denote by $\si(u)$ and $\tau(u)$ the first and the last letter of $u$, that is, $x_1$ and $x_n$, respectively. If $u$ and $v$ are two words with $\tau(u)=\si(v)$, we will denote by $u*v$ the word obtained by concatenating $u$ and $v$ without repeating $\si(u)=\tau(u)$ at their junction, that is, if $u=u_1x_n$ and $v=x_nv_1$, then $u*v=u_1x_nv_1$, for $x_n\in X$. 

A \emph{subsequence} $w$ of $u$ is any non-empty sequence of letters obtained from $u$ by deleting some (maybe none) of its letters. Thus two consecutive letters of $w$ may not be consecutive in $u$, but their relative position is maintained, that is, if $a$ is to the left of $b$ (not necessarily next to $b$) in $w$, then the same is true in $u$. A \emph{subword} $v$ of $u$ is a subsequence with `no holes', that is, $u=v_1vv_2$ for two possibly empty (with no letters) words $v_1$ and $v_2$. If $v_1$ is an empty word, we call $v$ a \emph{prefix} of $u$; and if $v_2$ is an empty word, we call $v$ a \emph{suffix} of $u$.

Let $R$ be a relation on $X^+$. A (semigroup) \emph{presentation} is a pair $\langle X,R\rangle$. The semigroup given by this presentation is the quotient semigroup $X^+/\rho$, where $\rho$ is the congruence generated by $R$. We will use in general the notation $[u]$ to refer the $\rho$-class $u\rho$, and we will write $u\approx v$ to indicate that $[u]=[v]$. The word problem for a presentation $\langle X,R\rangle$ is the question of knowing if $[u]=[v]$ for any $u,v\in X^+$.

Next we recall the Green's relations on a semigroup $S$. Since all our semigroups will be regular, we introduce these relations using their simplified version for regular semigroup. So, let $S$ be a regular semigroup. The quasi-orders $\leq_{\Lc}$,  $\leq_{\Rc}$ and $ \leq_{\Jc}$ on $S$ are defined as follows:
$$s\leq_{\Lc} t\Leftrightarrow Ss\subseteq St,\quad s\leq_{\Rc} t\Leftrightarrow sS\subseteq tS,\quad s\leq_{\Jc} t\Leftrightarrow SsS\subseteq StS.$$
Let $\geq_{\Lc}$, $\geq_{\Rc}$ and $\geq_{\Jc}$ be the corresponding dual relations. The first three Green's relations are
$$\Lc\,=\,\leq_{\Lc}\cap\geq_{\Lc},\quad \Rc\,=\,\leq_{\Rc}\cap\geq_{\Rc}\quad\mbox{ and }\quad \Jc\,=\,\leq_{\Jc}\cap\geq_{\Jc}\,.$$
The last two are $\Hc\,=\,\Lc\cap\Rc$ and $\Dc\,=\,\Lc\vee\Rc\,$. We will denote by $\Kcc_a$ the $\Kc$-class of the element $a\in S$, for $\Kc\in\{\Hc,\Lc,\Rc,\Dc,\Jc\}$. 

The \emph{natural partial order} $\leq$ on a regular semigroup $S$ is defined as follows:
$$s\leq t\quad\Leftrightarrow\quad  s=et=tf \;\mbox{ for some } e,f\in E(S)\,.$$
In this definition, we can also assume that $e\Rc s\Lc f$. This assumption does not change the relation $\leq\,$. Note also that $\leq\;\subseteq\;\leq_{\Rc}\cap\leq_{\Lc}$ but this inclusion is strict usually.

Another concept that will be very important to us is the notion of sandwich set. The \emph{sandwich set} of two idempotents $e$ and $f$ is the set
$$S(e,f)=\{g\in E(S)\,|\; fg=g=ge\; \mbox{ and }\; egf=ef\}\,.$$
This set is composed precisely by the inverses $g$ of $ef$ that satisfy $fg=g=ge$. If the semigroup $S$ is regular, the sets $S(e,f)$ are always nonempty. In general, the sandwich sets $S(e,f)$ and $S(f,e)$ are distinct sets, that is, we cannot permute the idempotents. However, $S(e,f)=S(e_1,f_1)$ for any $e,e_1,f,f_1\in E(S)$ such that $e\Lc e_1$ and $f\Rc f_1$. Hence, if $S$ is regular, we can extend the concept of sandwich set to non-idempotent elements: $S(a,b)=S(a'a,bb')$ for $a,b\in S$, $a'\in V(a)$ and $b'\in V(b)$.

\section{The presentation $\langle \G(X),\R(X)\rangle$}\label{sec4}

In this section we construct a presentation $\langle \G(X),\R(X)\rangle$. Due to its complexity, we divide it into four subsections to make it easy to follow. As we will see, most elements from $\G(X)$ will be 5-tuples. So, we need some special notation for dealing with these tuples. Also, all the definitions in this paper are related to a specific set $X$. To turn the notation less heavier we will omit the reference to the set $X$ in general. For example, instead of writing $\langle \G(X),\R(X)\rangle$, we will write only $\pres$. This will cause no ambiguity normally, but when it may cause, we will identify the sets $X$ in those cases. So, the reader should be always aware of this fact.

\subsection{Anchors and 5-tuples}

Given a nonempty set $X$, let
$$X'=\{x'\,|\;x\in X\}$$
be a disjoint copy of $X$. Let $\oX=X\cup X'$ and $A=\oX\cup\{1\}$, where $1$ is a new symbol. The elements of $A$ will be called \emph{anchors}. We establish an involution $'$ on $A$ by setting $1'=1$ and $(x')'=x$ for all $x\in X$. Then $x'$ will be called the \emph{inverse} (\emph{anchor}) of $x$.

Let $g$ be a 5-tuple. We will use the following notation to refer to the entries of $g$:
$$g=(g^l,g^{la},g^c,g^{ra},g^r)\,.$$
In the 5-tuples we will work with, the elements $g^{la}$ and $g^{ra}$ will be always anchors. Thus, we will call them the \emph{left} and \emph{right anchors} of $g$, respectively. The elements $g^l$, $g^c$ and $g^r$ will be again 5-tuples or $1$. We call these three elements, the \emph{left}, \emph{middle} and \emph{right entries} of $g$, respectively. We extend this 5-tuple notation to $g^l$, $g^c$ and $g^r$ as expected: for example
$$g^l=(g^{l^2},g^{l^2a},g^{lc},g^{lra},g^{lr})\,.$$
We believe it is now obvious how to use this notation in general (note that $g^{l^2}$, $g^{lc}$ and $g^{lr}$ will be often 5-tuples again, and so on).

\subsection{The sets $\G_{i,e}(X)$}\label{sub41}

The set $\G$ of generators will be constituted by 5-tuples, beside anchors. We will divide the set of 5-tuples of $\G$ into two subsets, the sets $\Gd$ and $\Ge$. The set $\Gd$ will contain the 5-tuples with distinct left and right entries, while $\Ge$ will contain the 5-tuples with equal left and right entries. Both $\Gd$ and $\Ge$ will be defined recursively as the union of subsets. We need to begin with $\Ge$ since these 5-tuples are needed to define $\Gd$. Thus, in this subsection, we define $\Ge$.

Let $\G_{1,e}=\{g_{xx'}\,|\;x\in X\}$ where
$$g_{xx'}=(1,1,xx',x',1)\,.$$
We define now $\Gie$ recursively as follows:
$$\Gie=\big\{\big(g,(g^{la})',g^l,(g^{ra})',g\big),\,\big(g,(g^{ra})',g^l,(g^{la})',g\big)\,:\;g\in\G_{i-1,e}\big\}$$
for $i\geq 2$. Note that $g^l=g^r$ since $g\in\G_{i-1,e}$, and therefore it is irrelevant to write $g^l$ or $g^r$ for the middle entry of the 5-tuples of $\Gie$. Further, $g^l\in\G_{i-2,e}$ if one considers $\G_{0,e}=\{1\}$. Hence
$$\Gie\subseteq \G_{i-1,e}\times A\times\G_{i-2,e}\times A\times\G_{i-1,e}$$
for $i\geq 2$. We set $\Ge=\cup_{i\in\mathbb{N}}\Gie$. We reinforce that we are not including $\G_{0,e}$ in $\Ge$ and so $\Ge$ is a set of 5-tuples only.

Analyzing more carefully the tuples from $\Ge$, the left entries completely determine the middle and right entries, that is, $g^c=g^{l^2}$ and $g^r=g^l$ for any $g\in\Gie$ with $i\geq 2$. The set $\G_{2,e}$ is constituted by the tuples
$$g_{2,x}=(g_{xx'},1,1,x,g_{xx'})\quad\mbox{ and }\quad g_{2,x'}=(g_{xx'},x,1,1,g_{xx'})\;,$$
for $x\in X$. The anchors are also determined by the index $i$:
$$\{g^{la},g^{ra}\}=\left\{\begin{array}{ll}
\{1,x\} &\mbox{ if } i \mbox{ even} \\ [.2cm]
\{1,x'\} &\mbox{ if } i \mbox{ odd\,.}
\end{array}\right.$$
The only alternative that exists is with their order: $g^{la}=1$ and $g^{ra}\neq 1$, or $g^{la}\neq 1$ and $g^{ra}=1$. Hence $\Gie$ has twice the number of elements of $\G_{i-1,e}\,$, and so $|\Gie|=2^{i-1}|X|$.

Note that both $g^{l^2}$ and $g^{lr}$ are equal to $g^c$ for each $g\in\Gie$ with $i\geq 2$. However, we need to attribute a `side' to $g^c$ inside $g^l$ with respect to $g$. Since only one of the elements $g^{l^2a}$ and $g^{lra}$ is equal to  $(g^{la})'$, we define
$$l_a=\left\{\begin{array}{ll}
	l^2 & \mbox{ if } g^{l^2a}=(g^{la})' \\ [.2cm]
	lr & \mbox{ if } g^{lra}=(g^{la})'\;.
\end{array}\right.$$
Hence $g^{l_a}=g^c$ and $g^{l_aa}=(g^{la})'$. But, although $g^{l_a}=g^c$, their meaning is different. While $g^c$ represents the middle entry of $g$, $g^{l_a}$ represents either the left or the right entry of $g^l$, that is, the `side' of $g^c$ inside $g^l$ with respect to $g$. When we use the notation $g^s$, especially $g^{l_a}$, we are not only considering its value but also its place inside a tuple.

In a similar manner we use the notation $r_a$. For each $g\in\Gie$ with $i\geq 2$, 
$$r_a=\left\{\begin{array}{ll}
	r^2 & \mbox{ if } g^{r^2a}=(g^{ra})' \\ [.2cm]
	rl & \mbox{ if } g^{rla}=(g^{ra})'\;.
\end{array}\right.$$
Then $g^{r_a}=g^c$ and $g^{r_aa}=(g^{ra})'$.

\subsection{The sets $\G_{i,d}(X)$}\label{sub42}

The next step is to introduce the 5-tuples of $\G$ with distinct left and right entries, that is, the elements of $\Gd$. We define $\Gd=\cup_{n\in\mathbb{N}}\Gid$ where the sets $\Gid$ are introduced recursively. However, to describe $\Gid$ we need not only the sets $\G_{i-1,d}$ but also $\G_{i-1,e}$. So let
$$\Gi=\Gid\cup\Gie$$
for $i\geq 1$ and set $\G_0=\{1\}$. 

The first set $\G_{1,d}$ is considered empty; whence $\G_1=\G_{1,e}$. Let $i\geq 2$ and assume that $\G_{i-1}$ and $\G_{i-2}$ are already defined and that the elements of $\G_{i-1}$ are 5-tuples. Set
$$\ol{\G}_i=\G_{i-1}\times A\times\G_{i-2}\times A\times\G_{i-1}\,.$$
Then $\Gid$ is the subset of $\ol{\G}_i$ constituted by the 5-tuples $g$ that satisfy
\begin{itemize}
	\item[$(i)$] $g^l\neq g^r$;
	\item[$(ii)$] $(g^c,g^{la})=(g^{ls},(g^{lsa})')$ for some $s\in\{l,r\}\,$; 
	\item[$(iii)$] $(g^c,g^{ra})=(g^{rt},(g^{rta})')$ for some $t\in\{l,r\}\,$. 
\end{itemize}
In particular, note that $g^c$ is always the left or right entry of $g^l$ and also the left or right entry of $g^r$.

The computation of the cardinality of the sets $\Gid$ is much harder than the case $\Gie\,$. Since this is not important for this paper, we will not do it here. We limit ourselves to describe the elements of $\G_{2,d}\,$. This list may turn the definition of these sets clearer. Its cardinality is $4|X|(|X|-1)$ since there are 4 tuples with left entry $g_{xx'}$ and right entry $g_{yy'}$ for each pair $(x,y)\in X\times X$ with $x\neq y$, namely
$$\begin{array}{ll}\label{G3d}
(g_{xx'},1,1,1,g_{yy'})\;,\qquad &(g_{xx'},x,1,1,g_{yy'})\;,\\ [.2cm]
(g_{xx'},1,1,y,g_{yy'})\;, \qquad &(g_{xx'},x,1,y,g_{yy'})\;.
\end{array}$$
As can be observe above, in the tuples $g\in\Gid$ we may have $g^{la}=g^{ra}$, something that does not occur in the tuples from $\Ge$. Nevertheless, we continue to know that $g^{la},g^{ra}\in\{1,x\}$ if $i$ even, and $g^{la},g^{ra}\in\{1,x'\}$ if $i$ odd, for some $x\in X$.

Let $g\in\Gid$. Since $g^l$ belongs to either $\G_{i-1,d}$ or $\G_{i-1,e}$, the definition of the `side' of $g^c$ inside $g^l$ with respect to $g$ is more complex. We define $l_a$ as  follows: $l_a=ls$ for $s\in\{l,r\}$ such that
$$g^{lsa}=(g^{la})' \;\mbox{ if } g^l\in\G_{i-1,e}\qquad \mbox{ or }\qquad  g^{ls}=g^c \;\mbox{ if } g^l\in\G_{i-1,d}\,.$$
Basically, if $g^l\in\G_{i-1,e}$, and so $g^{l^2}=g^c=g^{lr}$, we define $l_a$ as in the previous subsection; if $g^l\in\G_{i-1,d}$, the `side' of $g^c$ inside $g^l$ is naturally defined by the fact that only one of the entries $g^{l^2}$ or $g^{lr}$ of $g^l$ is equal to $g^c$. 

We define $r_a$ similarly. Note that
$$g^{l_aa}=(g^{la})'\quad\mbox{ and }\quad g^{r_aa}=(g^{ra})'$$
again for the elements of $\Gid$. Thus, the previous equalities hold for all cases where $l_a$ is defined (that is, except for the elements of $\G_{1,e}$). We alert that these equalities will be used often in this paper without further comments.

Let $\Gf=\cup_{i\in\mathbb{N}}\Gi=\Gd\cup\Ge$. The \emph{height} $\up(g)$ of $g\in\Gf$ is the index $i$ such that $g\in\Gi$. Set $\G=\Gf\cup A$ and $\G'=\Gf\cup\{1\}$. Consider $1$ to have height $0$.

\subsection{The relation $\R(X)$}\label{sub43}

Let $g\in\Gf$ of height $\up(g)=i$ at least $2$. Let $g^L$ and $g^R$ be the triplets
$$g^L=(g^{la})'g^lg^{la}\qquad\mbox{ and }\qquad g^R=(g^{ra})'g^rg^{ra}\,.$$
We define $\R$ as the union $\rho_e\cup\rho_s$ where
$$\rho_e=\{(xx'x,x),(x'xx'\!,x'),(g_{xx'},xx'),(1g,g),(g1,g),(g^2,g)|\,x\!\in\! X,\,g\!\in\!\G'\}$$
and
$$\rho_s=\left\{(g^cg^Lg,g),(gg^Rg^c,g),(g^Rgg^L,g^Rg^cg^L)|\, g\in \G^5 \mbox{ with }\up(g)\geq 2\right\}.$$
We denote by $\fox$ the semigroup given by the presentation $\pres$, that is, $\fox=\Gp/\rho$ where $\rho$ is the smallest congruence containing $\R$.

Note that the two first pairs of $\rho_e$ indicate that $[x']$ is an inverse of $[x]$, the third that $[g_{xx'}]=[xx']$, the forth and fifth that $[1]$ is an identity element for $\fox$, and finally the last that the elements $[g]$ are idempotents of $\fox$. Thus $\fox$ is an idempotent generated monoid with identity element $[1]$. This is the meaning of $\rho_e$. 

The next result explains the meaning of $\rho_s$ (Lemma \ref{prodelem}.$(v)$) and it is the generalization of \cite[Lemma 3.1]{LO23}. The proof is the same except for $(v)$ where the small difference comes from the fact the definition of $\rho_s$ here is not exactly the same as in \cite{LO23} (though equivalent). Hence, we do next only the proof of $(v)$. 

\begin{lem}\label{prodelem}
Let $g\in \G$ with $\up(g)>1$. Then:
\begin{itemize}
\item[$(i)$] $g\ap g^cg\ap gg^c\ap gg^Lg\ap gg^Rg$.
\item[$(ii)$] $[gg^L],\,[gg^R],\,[g^Lg],\,[g^Rg]\in E(\fox)$.
\item[$(iii)$] $g^Lg^cg^L\ap g^L$ and $g^Rg^cg^R\ap g^R$. 
\item[$(iv)$] $[g^cg^L],\,[g^cg^R],\,[g^Lg^c],\,[g^Rg^c]\in E(\fox)$.
\item[$(v)$] $[g]\in S([g^Rg^c],[g^cg^L])$.
\end{itemize}
\end{lem}

\begin{proof}
$(v)$. Note that we only need to verify that $g^Rg^cgg^cg^L\ap (g^Rg^c)(g^cg^L)\,$:
$$g^Rg^cgg^cg^L\ap g^Rgg^L\ap g^Rg^cg^L\ap (g^Rg^c)(g^cg^L)\,;$$
whence $[g]\in S([g^Rg^c],[g^cg^L])$.
\end{proof}

We have defined the presentation $\pres$ and explained the purpose of the pairs in the relation $\R$. In the next section we solve the word problem for this presentation.

\section{The word problem for $\langle\G(X),\R(X)\rangle$}
 
In order to facilitate the understanding and the notation that follows, we will use $a$ and $b$ to refer to anchors, $g$ and $h$ to refer to letters of $\G'$, and $u$ and $v$ to refer to words from $\G^+$ during this paper. Of course, we may include also indices in these letters.

\subsection{Special words}\label{sub51}

A triplet $g_1ag_2$ with $g_1,g_2\in\G'$ and $a\in A$ is \emph{left anchored} if $$(g_1,a)=(g_2^s,g_2^{sa})\;\mbox{ for some }\;s\in\{l,r\}\,,$$ 
and \emph{right anchored} if
$$(g_2,a)=(g_1^s,(g_1^{sa})')\;\mbox{ for some }\;s\in\{l,r\}\,.$$
In the former case we say that $g_1$ is \emph{anchored} to $g_2$ and in the latter that $g_2$ is \emph{anchored} to $g_1$. We reinforce here that these notions are not exactly self dual: for the left anchored case, $a$ is an anchor of $g_2$; while for the right anchored case, $a$ is the inverse of an anchor of $g_1$. In general, a triplet $g_1ag_2$ is \emph{anchored} if it is left or right anchored.

A word $u=g_0a_1g_1\cdots a_ng_n\in \Gp$ with $n\geq 1$ is a \emph{landscape} if each triplet $g_{i-1}a_ig_i$ is anchored for $1\leq i\leq n$. We consider also the letters from $\G'$ as landscapes. If $n\geq 1$, the subsequences $a_1\cdots a_n$ and $g_0g_1\cdots g_n$ are called the \emph{anchors subsequence} and the \emph{letters subsequence} of $u$. We denote by $\lx$ the set of all landscapes of $\Gp$.

\begin{figure}[ht]
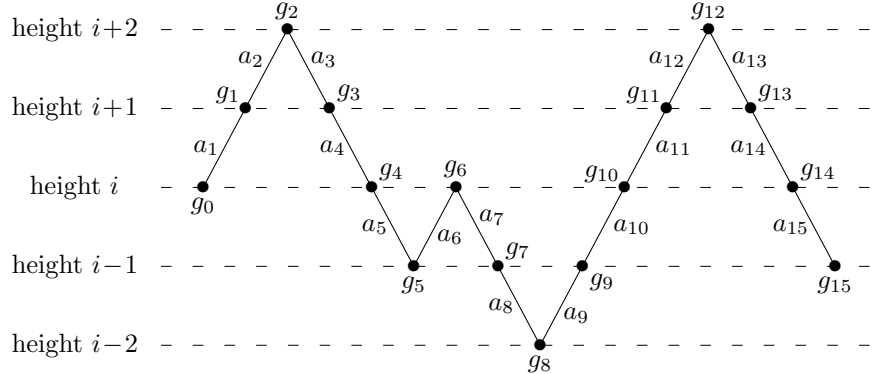
  
$$\tikz[scale=.7]{
\coordinate (1) at (1.8,3);
\coordinate (2) at (2.6,4.5);
\coordinate (3) at (3.4,6);
\coordinate (4) at (4.2,4.5);
\coordinate (5) at (5,3);
\coordinate (6) at (5.8,1.5);
\coordinate (7) at (6.6,3);
\coordinate (8) at (7.4,1.5);
\coordinate (9) at (8.2,0);
\coordinate (10) at (9,1.5);
\coordinate (11) at (9.8,3);
\coordinate (12) at (10.6,4.5);
\coordinate (13) at (11.4,6);
\coordinate (14) at (12.2,4.5);
\coordinate (15) at (13,3);
\coordinate (16) at (13.8,1.5);
\coordinate (17) at (1,0);
\coordinate (18) at (1,1.5);
\coordinate (19) at (1,3);
\coordinate (20) at (1,4.5);
\coordinate (21) at (1,6);
\coordinate (22) at (14.5,0);
\coordinate (23) at (14.5,1.5);
\coordinate (24) at (14.5,3);
\coordinate (25) at (14.5,4.5);
\coordinate (26) at (14.5,6);
\draw (1) node {$\bullet$} node [below] {\small $g_0$};
\draw (2) node {$\bullet$} node [above left=-2pt] {\small $g_1$};
\draw (3) node {$\bullet$} node [above] {\small $g_2$};
\draw (4) node {$\bullet$} node [above right=-1pt] {\small $g_3$};
\draw (5) node {$\bullet$} node [above right=-1pt] {\small $g_4$};
\draw (6) node {$\bullet$} node [below] {\small $g_5$};
\draw (7) node {$\bullet$} node [above] {\small $g_6$};
\draw (8) node {$\bullet$} node [above right=-1pt] {\small $g_7$};		
\draw (9) node {$\bullet$} node [below] {\small $g_8$};		
\draw (10) node {$\bullet$} node [below right=-1pt] {\small $g_9$};
\draw (11) node {$\bullet$} node [above left=-2pt] {\small $g_{10}$};
\draw (12) node {$\bullet$} node [above left=-2pt] {\small $g_{11}$};
\draw (13) node {$\bullet$} node [above] {\small $g_{12}$};	
\draw (14) node {$\bullet$} node [above right=-1pt] {\small $g_{13}$};
\draw (15) node {$\bullet$} node [above right=-1pt] {\small $g_{14}$};		
\draw (16) node {$\bullet$} node [below] {\small $g_{15}$};						
\draw (17) node [left=5pt] {\small height $i\!-\!2$};
\draw (18) node [left=5pt] {\small height $i\!-\!1$};
\draw (19) node [left=12pt] {\small height $i$};
\draw (20) node [left=5pt] {\small height $i\!+\!1$};
\draw (21) node [left=5pt] {\small height $i\!+\!2$};		
\draw[very thin, dashed, dash pattern=on 4pt off 8pt] (17)--(22);
\draw[very thin, dashed, dash pattern=on 4pt off 8pt] (18)--(23);
\draw[very thin, dashed, dash pattern=on 4pt off 8pt] (19)--(24);
\draw[very thin, dashed, dash pattern=on 4pt off 8pt] (20)--(25);
\draw[very thin, dashed, dash pattern=on 4pt off 8pt] (21)--(26);
\draw (1) to node[left=-2pt] {\small$a_1$} (2);
\draw (2) to node[above left=-3pt] {\small$a_2$} (3);
\draw (3) to node[above right=-3pt] {\small$a_3$} (4);
\draw (4) to node[left=-2pt] {\small$a_4$} (5);
\draw (5) to node[left=-2pt] {\small$a_5$} (6);
\draw (6) to node[below right=-3pt] {\small$a_6$} (7);
\draw (7) to node[above right=-3pt] {\small$a_7$} (8);
\draw (8) to node[left=-2pt] {\small$a_8$} (9);
\draw (9) to node[below right=-3pt] {\small$a_9$} (10);
\draw (10) to node[right] {\small$a_{10}$} (11);
\draw (11) to node[right] {\small$a_{11}$} (12);
\draw (12) to node[above left=-3pt] {\small$a_{12}$} (13);
\draw (13) to node[above right=-3pt] {\small$a_{13}$} (14);
\draw (14) to node[left=-2pt] {\small$a_{14}$} (15);	
\draw (15) to node[left=-2pt] {\small$a_{15}$} (16);				
}$$
\caption{Illustration of a landscape in a line graph.}\label{fig1}
\end{figure}

We will represent the landscapes $u=g_0a_1g_1\cdots a_ng_n$ in a line graph with vertices $\{g_0,g_1,\cdots,g_n\}$ and edges $\{a_1,\cdots,a_n\}$. As should be expected, the edge $a_i$ connects the vertices $g_{i-1}$ and $g_i$. We add also information about the height of the letters by drawing vertices with the same height in the same imaginary horizontal line. Further, higher letters should appear in higher imaginary horizontal lines. The anchored triplet $g_{i-1}a_ig_i$ is left anchored if $g_i$ is higher than $g_{i-1}$, and right anchored if $g_{i-1}$ is higher than $g_i$. We include here the figure used in \cite{LO23} to illustrate a landscape (Figure \ref{fig1}). 

The reader may have notice that we are using $g_{xx'}$ to represent $xx'$ (that is, $[g_{xx'}]=[xx']$) and may found it awkward and probably unnecessary. However, the definition of landscape gives us a reason for this option. When dealing with landscapes, we want to look at $xx'$ not as the concatenation of two anchors but instead as a letter from $\G'$. Thus writing $g_{xx'}$ instead of $xx'$ helps making this distinction.

\begin{figure}[ht]
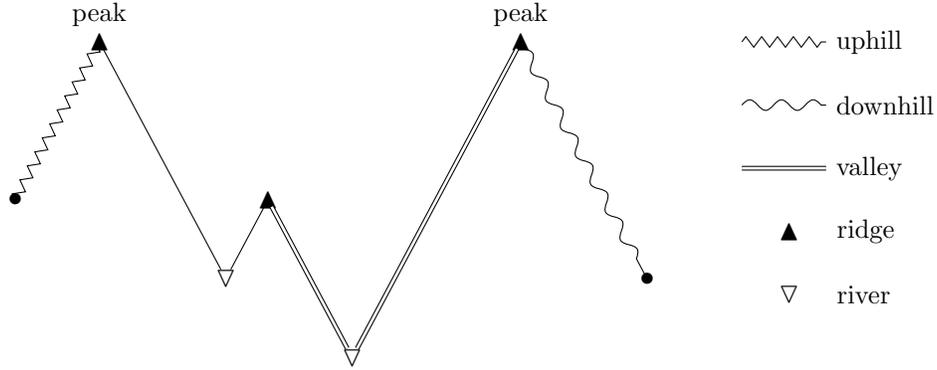
 
$$\tikz[scale=.7]{
\coordinate (1) at (1.8,3);
\coordinate (2) at (3.4,6);
\coordinate (3) at (5.8,1.5);
\coordinate (4) at (6.6,3);
\coordinate (5) at (8.2,0);
\coordinate (6) at (11.4,6);
\coordinate (7) at (13.8,1.5);
\draw (1) node {$\bullet$};
\draw (2) node {$\blacktriangle$} node[above=2pt] {\small peak};
\draw (3) node {$\triangledown$};
\draw (4) node {$\blacktriangle$};
\draw (5) node {$\triangledown$};	
\draw (6) node {$\blacktriangle$}	node[above=2pt] {\small peak};
\draw (7) node {$\bullet$};	
\draw[snake,segment amplitude=2pt,segment length=6pt,gap after snake=1pt] (1)--(2);
\draw[shorten >=4pt, shorten <=2pt] (2)--(3);
\draw[shorten >=2pt, shorten <=4pt] (3)--(4);
\draw[double,double distance=1pt,shorten >=4pt, shorten <=3pt] (4)--(5);
\draw[double,double distance=1pt,shorten >=3pt, shorten <=4pt] (5)--(6);
\draw[snake=coil,segment aspect=0,segment amplitude=2pt,segment length=12pt,gap before snake=3pt] (6)--(7);			
\draw[snake,segment amplitude=2pt,segment length=6pt] (15.6,6)--(17.2,6);
\draw (17.2,6) node[right] {\small uphill};
\draw[snake=coil,segment aspect=0,segment amplitude=2pt,segment length=12pt] (15.6,4.8)--(17.2,4.8);
\draw (17.2,4.8) node[right] {\small downhill};
\draw[double,double distance=1pt] (15.6,3.6)--(17.2,3.6);
\draw (17.2,3.6) node[right] {\small valley};
\draw (16.5,2.4) node {$\blacktriangle$} node[right=.5cm] {\small ridge};
\draw (16.5,1.2) node {$\triangledown$} node[right=.5cm] {\small river};
}$$
\caption{Illustration of some terminology.}\label{fig2}
\end{figure}

The graphical representation of landscapes is very useful now to introduce other terms such as uphills, downhills, valleys, rivers, ridges and peaks. We will use these terms both for words and for subwords of other words. We illustrate these notions for subwords in Figure \ref{fig2} (a figure also used in \cite{LO23}). We believe this figure is enough to explain these terms. Nevertheless we give next their formal definition.

A landscape $u=g_0a_1g_1\cdots a_ng_n$ is called an \emph{uphill} if $\up(g_i)=\up(g_{i-1})+1$ for all $1\leq i\leq n$ and a \emph{downhill} if $\up(g_{i})=\up(g_{i-1})-1$ for all $1\leq i\leq n$. A \emph{hill} is then a landscape falling into one of these two special cases. Let $\lo$, $\lop$ and $\lom$ denote the set of all hills, uphills and downhills, respectively. We denote also by $\ltp$ and $\ltm$ the sets of all landscapes composed, respectively, by an uphill followed by a downhill and by a downhill followed by an uphill. Set $\lt=\ltp\cup\ltm$. The words from $\ltm$ are called \emph{valleys}. The members of $\ltp$ have no special name except for a specific kind, the mountains. A nontrivial \emph{mountain} is a landscape $u$ from $\ltp$ with $\si(u)=1=\tau(u)$. The trivial mountain is the word $1$ and, of course, is the only mountain not in $\ltp$.

A letter $g_i$ with $0<i<n$ is a \emph{river} [\emph{ridge}] of the landscape $u$ if $\up(g_{i-1})=\up(g_{i+1})=\up(g_i)+1$ [$\up(g_{i-1})=\up(g_{i+1})=\up(g_i)-1$]. A \emph{peak} of $u$ is a ridge of highest height. Note that a landscape can have several distinct peaks but they all have the same height. If $u$ has a unique peak, we denote it by $\ka(u)$. We generalize the notion of height for letters of $\G'$ to all landscapes by defining the height of u, denoted by $\up(u)$, as the maximum between the height of its peaks, of its first letter and of its last letter.

A \emph{mountain range} is a landscape $u$ with $\si(u)=1=\tau(1)$. Thus the mountains are the mountain ranges with no rivers and the nontrivial mountains are the mountain ranges composed by an uphill followed by a downhill. If $u$ is a nontrivial mountain range, we denote by $\la_l(u)$ and $\la_r(u)$, respectively, its initial maximal uphill and its final maximal downhill. We call $\la_l(u)$ the \emph{left hill} of $u$ and $\la_r(u)$ the right hill of $u$. For completeness we set also $\la_l(1)=\la_r(1)=1$. The height of a nontrivial mountain range is clearly the height of its peaks. Hence, if $u$ is a nontrivial mountain, then $u=\la_l(u)*\la_r(u)$ and $\up(u)=\up(\ka(u))$. We denote by $\mrx$, $\mxm$ and $\mx$ the sets of all mountain ranges, mountains and nontrivial mountains, respectively. Then $\mx=\mrx\cap \ltp$. Finally, note that all mountain ranges have length $4k+1$ for some $k\in\mathbb{N}_0$.

If $u$ and $v$ are two landscapes with $\tau(u)=g=\si(v)$, the concatenation $uv$ is not a landscape but $u*v$ is: there is a double $g$ ($gg$) as a subword of $uv$. Note also that $[uv]=[u*v]$. Further, $u*v\in\mrx$ if $u,v\in\mrx$. 

Next, we introduce notation for special anchored triplets. We write $g^s\cdot g$ and $g\cdot g^s$ to refer to the anchored triplets $g^sg^{sa}g$ and $g(g^{sa})'g^s$, respectively, for $s\in\{l,r\}$. We alert that this `dot' notation can be used only when the value of $s$ is explicitly indicated. Otherwise, it may become ambiguous: if $g\in\G_e$ and $h$ is then simultaneously its left and right entry, then $h\cdot g$ could refer to either $hg^{la}g$ and $hg^{ra}g$. We write also $g^c\cdot g^s$ and $g^s\cdot g^c$ to refer to the anchored triplets 
$$g^cg^{s_aa}g^s=g^c(g^{sa})'g^s\quad\mbox{ and }\quad g^s(g^{s_aa})'g^c= g^sg^{sa}g^c\,,$$ 
respectively. Here again the indication of the sup-script $c$ and the specification of $s$ are absolutely necessary to avoid ambiguities.

\begin{lem}\label{l51a}
For each $g\in\G^5$ with $\up(g)>1$, $\;g^r\cdot g\cdot g^l\ap g^r\cdot g^c\cdot g^l$.
\end{lem}

\begin{proof}
We begin by showing that $g^lg^{la}g^c(g^{la})'g^l\ap g^l$. Note that
$$g^lg^{la}g^c(g^{la})'g^l=g^l(g^{l_aa})'g^cg^{l_aa}g^l=\left\{\begin{array}{ll}
g^l(g^l)^Lg^l & \mbox{ if } l_a=l^2 \\ [.2cm]
g^l(g^l)^Rg^l & \mbox{ if } l_a=lr\,.
\end{array}\right.$$
By Lemma \ref{prodelem}.$(i)$ applied to $g^l$ we now conclude that $g^lg^{la}g^c(g^{la})'g^l\ap g^l$. Similarly, we can also state that $g^rg^{ra}g^c(g^{ra})'g^r\ap g^r$. Finally,
$$\begin{array}{lllr}
g^r\cdot g\cdot g^l&=&g^rg^{ra}g(g^{la})'g^l & \\
&\ap&g^rg^{ra}g^c(g^{ra})'g^rg^{ra}g(g^{la})'g^lg^{la}g^c(g^{la})'g^l & \\
&=&g^rg^{ra}g^cg^Rgg^Lg^c(g^{la})'g^l & \\
&\ap&g^rg^{ra}g^cg^Rg^cg^Lg^c(g^{la})'g^l & (\mbox{by definition of }\rho_s)\\
&\ap&\cdots\;\;\ap\;\; g^rg^{ra}g^c(g^{la})'g^l=g^r\cdot g^c\cdot g^l\,.
\end{array}$$
We have proved this lemma.
\end{proof}

The words $g^c\cdot g^l\cdot g$ and $g\cdot g^r\cdot g^c$ are very special to us since they constitute the `building blocks' to construct an important mountain, the mountain $\be_1(g)$. But first note that
$$g^c\cdot g^l\cdot g=g^cg^Lg\ap g\quad\mbox{ and } \quad g\cdot g^r\cdot g^c=gg^Rg^c\ap g\,.$$
Let $g\in\G^5$ and define
$$\la_l(g)=\left\{\begin{array}{ll}
\underbracket{\,g^{c^n}\cdot g^{c^{n-1}l}\cdot\,}g^{c^{n-1}}\!\cdots \underbracket{\,g^{c^2}\cdot g^{cl}\cdot\,}\underbracket{\,g^c\cdot g^l\cdot\,}g & \mbox{ if } \up(g)=2n \\ [.3cm]
11\underbracket{\,g^{c^n}\cdot g^{c^{n-1}l}\cdot\,}g^{c^{n-1}}\!\cdots \underbracket{\,g^{c^2}\cdot g^{cl}\cdot\,}\underbracket{\,g^c\cdot g^l\cdot\,}g & \mbox{ if } \up(g)=2n+1\,.
\end{array}\right.$$
The under-brackets try to emphasize the `building blocks' $g^{c^i}\cdot g^{c^{i-1}l}\cdot g^{c^{i-1}}$ used in the construction of $\la_l(g)$. Then $\la_l(g)$ is an uphill from $1$ to $g$ such that $\la_l(g)\ap g$ by recursion. 

\begin{figure}[ht]
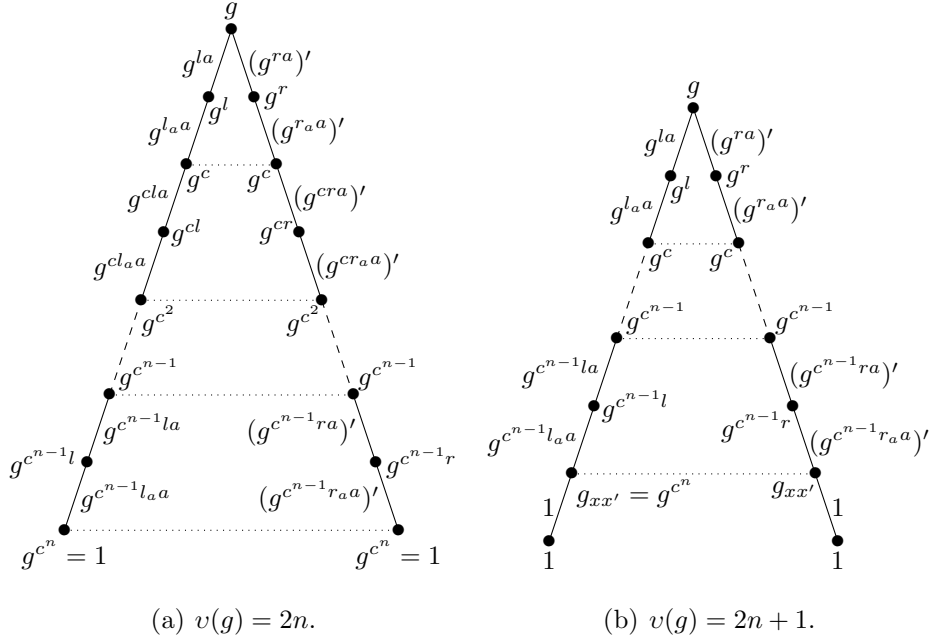

\begin{subfigure}[b]{.49\textwidth}
\centering
$$\tikz[scale=.6]{
\coordinate (1) at (0,0);
\coordinate (2) at (.5,1.5);
\coordinate (3) at (1.2,3.6);
\coordinate (4) at (1.7,5.1);
\coordinate (5) at (2.2,6.6);
\coordinate (6) at (2.7,8.1);
\coordinate (7) at (3.2,9.6);
\coordinate (8) at (3.7,8.1);
\coordinate (9) at (4.2,6.6);
\coordinate (10) at (4.7,5.1);
\coordinate (11) at (5.2,3.6);
\coordinate (12) at (5.9,1.5);
\coordinate (13) at (6.4,0);
\coordinate (14) at (-0.5,-1.5);
\coordinate (15) at (6.9,-1.5);		
\draw (14) node {$\bullet$} node [below] {\small $g^{c^n}=1$};		
\draw (1) node {$\bullet$} node [left] {\small $g^{c^{n-1}l}$};
\draw (2) node {$\bullet$} node [above right=-2pt] {\small $\,g^{c^{n-1}}$};		
\draw (3) node {$\bullet$} node [below right=-2pt] {\small $\!g^{c^2}$};
\draw (4) node {$\bullet$} node [right=-1pt] {\small $g^{cl}$};
\draw (5) node {$\bullet$} node [below right=-2pt] {\small $\!g^c$};		
\draw (6) node {$\bullet$} node [below right=-4pt] {\small $g^l$};		
\draw (7) node {$\bullet$} node [above] {\small $g$};
\draw (8) node {$\bullet$} node [right] {\small $g^r$};		
\draw (9) node {$\bullet$} node [below left=-2pt] {\small $g^c$};
\draw (10) node {$\bullet$} node [left=-2pt] {\small $g^{cr}$};
\draw (11) node {$\bullet$} node [below left=-2pt] {\small $g^{c^2}\!$};		
\draw (12) node {$\bullet$} node [above right=-2pt] {\small $g^{c^{n-1}}$};			
\draw (13) node {$\bullet$} node [right] {\small $g^{c^{n-1}r}$};
\draw (15) node {$\bullet$} node [below] {\small $g^{c^n}=1$};		
\draw (14) to node[right=-2pt] {\small$g^{c^{n-1}l_aa}$} (1);
\draw (1) to node[right=-2pt] {\small$g^{c^{n-1}la}$} (2);		
\draw[dashed] (2)--(3);
\draw (3) to node[left=-2pt] {\small$g^{cl_aa}$} (4);
\draw (4) to node[left=-2pt] {\small$g^{cla}$} (5);
\draw (5) to node[left=-2pt] {\small$g^{l_aa}$} (6);
\draw (6) to node[left=-2pt] {\small$g^{la}$} (7);
\draw (7) to node[right=-2pt] {\small$(g^{ra})'$} (8);
\draw (8) to node[right=-2pt] {\small$(g^{r_aa})'$} (9);
\draw (9) to node[right=-2pt] {\small$(g^{cra})'$} (10);
\draw (10) to node[right=-2pt] {\small$(g^{cr_aa})'$} (11);
\draw[dashed] (11)--(12);
\draw (12) to node[left=-1pt] {\small$(g^{c^{n-1}ra})'$} (13);
\draw (13) to node[left=-1pt] {\small$(g^{c^{n-1}r_aa})'$} (15);				
\draw[dotted] (14)--(15);
\draw[dotted] (3)--(11);
\draw[dotted] (5)--(9);		
\draw[dotted] (2)--(12);				
}$$	
\caption{$\up(g)=2n$.}
\end{subfigure}
\begin{subfigure}[b]{.49\textwidth}
\centering
$$\tikz[scale=.6]{
\coordinate (1) at (0,0);			
\coordinate (2) at (0.5,1.5);
\coordinate (3) at (1,3);
\coordinate (4) at (1.5,4.5);
\coordinate (5) at (2.2,6.6);
\coordinate (6) at (2.7,8.1);
\coordinate (7) at (3.2,9.6);
\coordinate (8) at (3.7,8.1);
\coordinate (9) at (4.2,6.6);
\coordinate (10) at (4.9,4.5);
\coordinate (11) at (5.4,3);
\coordinate (12) at (5.9,1.5);
\coordinate (13) at (6.4,0);		
\draw (1) node {$\bullet$} node [below] {\small $1$};
\draw (2) node {$\bullet$} node [below right=-2pt] {\small $g_{xx'}=g^{c^n}$};		
\draw (3) node {$\bullet$} node [right=-1pt] {\small $g^{c^{n-1}l}$};
\draw (4) node {$\bullet$} node [above right=-2pt] {\small $\,g^{c^{n-1}}$};		
\draw (5) node {$\bullet$} node [below right=-2pt] {\small $\!g^c$};		
\draw (6) node {$\bullet$} node [below right=-4pt] {\small $g^l$};		
\draw (7) node {$\bullet$} node [above] {\small $g$};
\draw (8) node {$\bullet$} node [right] {\small $g^r$};		
\draw (9) node {$\bullet$} node [below left=-2pt] {\small $g^c$};
\draw (10) node {$\bullet$} node [above right=-2pt] {\small $g^{c^{n-1}}$};			
\draw (11) node {$\bullet$} node [below left=-2pt] {\small $g^{c^{n-1}r}\!$};
\draw (12) node {$\bullet$} node [below left=-1pt] {\small $g_{xx'}\!\!$};
\draw (13) node {$\bullet$} node [below] {\small $1$};		
\draw (1) to node[left=-2pt] {\small$1$} (2);
\draw (2) to node[left=-2pt] {\small$g^{c^{n-1}l_aa}$} (3);
\draw (3) to node[left=-2pt] {\small$g^{c^{n-1}la}$} (4);		
\draw[dashed] (4)--(5);
\draw (5) to node[left=-2pt] {\small$g^{l_aa}$} (6);
\draw (6) to node[left=-2pt] {\small$g^{la}$} (7);
\draw (7) to node[right=-2pt] {\small$(g^{ra})'$} (8);
\draw (8) to node[right=-2pt] {\small$(g^{r_aa})'$} (9);
\draw[dashed] (9)--(10);
\draw (10) to node[right=-1pt] {\small$(g^{c^{n-1}ra})'$} (11);
\draw (11) to node[right=-2pt] {\small$(g^{c^{n-1}r_aa})'$} (12);	
\draw (12) to node[right=-2pt] {\small$1$} (13);			
\draw[dotted] (2)--(12);
\draw[dotted] (5)--(9);		
\draw[dotted] (4)--(10);				
}$$	
\caption{$\up(g)=2n+1$.}
\end{subfigure}
\caption{Illustration of the mountain $\be_1(g)$.}\label{fig3}
\end{figure}

We define similarly the downhill $\la_r(g)$ from $g$ to $1$ as follows:
$$\la_l(g)=\left\{\begin{array}{ll}
g\underbracket{\,\cdot g^r\cdot g^c\,}\underbracket{\,\cdot g^{cr}\cdot g^{c^2}\,}\!\cdots g^{c^{n-1}}\underbracket{\,\cdot g^{c^{n-1}r}\cdot g^{c^n}} & \mbox{ if } \up(g)=2n \\ [.3cm]
g\underbracket{\,\cdot g^r\cdot g^c\,}\underbracket{\,\cdot g^{cr}\cdot g^{c^2}\,}\!\cdots g^{c^{n-1}}\underbracket{\,\cdot g^{c^{n-1}r}\cdot g^{c^n}}11 & \mbox{ if } \up(g)=2n+1\,.
\end{array}\right.$$
Then $\la_r(g)\ap g$ again by recursion. Define $\be_1(g)=\la_l(g)*\la_r(g)$, a mountain with peak $g$ such that $\be_1(g)\ap g$. To help understanding the construction of $\be_1(g)$, we illustrate it in Figure \ref{fig3}, again a figure already used in \cite{LO23} for the same purpose. Note also that $\be_1(g_{xx'})=11g_{xx'}11$ for any $x\in X$. For completeness we define also $\be_1(1)=1$.

It may look strange the appearance of $11$ as subwords of $\be_1(g)$ since we know that $[1]$ is the identity element of $\fox$, and so the double $1$ could be omitted or replaced by a single $1$. The reason for writhing $11$ is because one of them is considered as a letter from $\G'$ and the other an anchor. This is crucial to guarantee that we have a mountain range in the way we have defined it. Any alternative definition of landscape or mountain range that would allow us to avoid writing $11$, would make us treat $1$ (or the letters of height $1$ from $\G'$) differently from the other letters of $\G'$. So, writing these double $1$ is the technical way to go in order to treat all letters from $\G'$ equally, and avoid having to analyze the case of the letter $1$ apart from the others.  

We need to extend $\be_1$ to all words from $\G^+$. We begin by extending $\be_1$ to the set $\oX$. Note that $11g_{xx'}x1$ and $1x'g_{xx'}11$ are mountains, and 
$$11g_{xx'}x1\ap x\quad\mbox{ and }\quad 1x'g_{xx'}11\ap x'\,,$$
for each $x\in X$. Naturally we define 
$$\be_1(x)=11g_{xx'}x1\quad\mbox{ and }\quad\be_1(x')=1x'g_{xx'}11\,.$$ 
We extend now $\be_1$ naturally to all words from $\G^+$ by setting
$$\be_1(e_1e_1\cdots e_n)=\be_1(e_1)*\be_1(e_2)*\cdots *\be_1(e_n)\,,$$ 
where $e_i\in\G$ for $1\leq i\leq n$. The words $\be_1(u)$ for $u\in\G^+$ are clearly mountain ranges (not mountains unless $n=1$) since all $\be_1(e_i)$ are mountains. The next result is now obvious from the observations already made:

\begin{lem}\label{l52a}
For each $u\in\G^+$, $\be_1(u)\ap u$ and $\be_1(u)$ is a mountain range.
\end{lem}

Let $g\in\G_n$ and let $u=g_0a_1g_1\cdots g_{n-1}a_ng_n$ be an uphill from $1$ to $g$. Then $g_0=1$ and $g_n=g$. Let also $i$ be the largest integer such that $g_i\in\G_e$. Clearly $i\geq 1$ since $g_1\in\G_e$ always. Note that one of the anchors of $g_i$ is $1$ and the other belongs to $\oX$. Consider $a\in X$ such that either $a$ or $a'$ is the anchor of $g_i$ distinct from $1$. By definition of the tuples from $\G^5$, $g_j\in\G_e$ if $j\leq i$ and $g_j\in\G_d$ if $j>i$. If $j>i$ and $g_j$ is known, then there are exactly two distinct options for $g_{j-1}$, either $g_j^l$ or $g_j^r$. Moreover, the anchor $a_i$ is fixed once $g_{j-1}$ is chosen. If $j\leq i$, then $g_{j-1}$ is determined by $g_j$ since $g_j^l=g_j^r$, but we have two options for the anchor $a_i$, one of them is $1$, the other is either $a$ or $a'$ depending on $j$ being even or odd, respectively. Therefore, we can conclude there are $2^n$ uphills from $1$ to $g$.

We could give a similar description for the downhills from $g$ to $1$, but we will leave it for the reader. We will take instead this chance to introduce another important concept, the concept of reverse of a landscape. If $g_1ag_2$ is an anchored triplet, then $g_2ag_1$ is not anchored if $a\neq 1$. However $g_2a'g_1$ is always anchored. Thus, if $u=g_0a_1g_1\cdots g_{n-1}a_ng_n$ is a landscape, its \emph{reverse}
$$\cev{u}=g_na_n'g_{n-1}\cdots g_1a_1'g_0$$
is another landscape; and the reverse of $\cev{u}$ is $u$ again. This `reverse operation' gives us a natural bijection between the sets of uphills from $1$ to $g$ and of downhills from $g$ to $1$. So, there are also $2^n$ downhills from $g$ to $1$. If we now consider the mountains with peak $g$, there must exist exactly $2^{2n}$ such mountains by the previous discussion. We have shown the next result:

\begin{lem}\label{number}
Let $g\in\G_n$. There are $2^{2n}$ mountains with peak $g$, $2^n$ uphills from $1$ to $g$, and $2^n$ downhills from $g$ to $1$.
\end{lem}

In the next subsection we introduce a procedure that will transform each mountain range into a mountain.

\subsection{Uplifting of rivers}\label{sub52}

Let $u=g_0a_1g_1\cdots a_ng_n$ be a landscape and assume $g_i$ is a river of $u$. Then $\up(g_{i-1})=\up(g_{i+1})=\up(g_i)+1$, $g_{i-1}a_ig_i$ is right anchored triplet and $g_ia_{i+1}g_{i+1}$ is a left anchored triplet. Hence,
$$(g_i,a_i)=\big(g_{i-1}^s,(g_{i-1}^{sa})'\big)\quad\mbox{ and }\quad (g_i,a_{i+1})= (g_{i+1}^t,g_{i+1}^{ta})\,,$$
for some $s,t\in\{l,r\}$. If $(g_{i-1},a_i)=(g_{i+1},a_{i+1}')$, then $g_{i-1}a_ig_ia_{i+1}g_{i+1}$ is either $g_{i-1}g_{i-1}^Lg_{i-1}$ or $g_{i-1}g_{i-1}^Rg_{i-1}$; whence $g_{i-1}a_ig_ia_{i+1}g_{i+1}\ap g_{i-1}$ by Lemma \ref{prodelem}.$(i)$. If $(g_{i-1},a_i)\neq (g_{i+1},a_{i+1}')$, then
$$h_i=(g_{i+1},a_{i+1}',g_i,a_i,g_{i-1})\in\G^5\,.$$
More precisely, $h_i\in\G_d$ if $g_{i-1}\neq g_{i+1}$ and $h_i\in\G_e$ if $g_{i-1}= g_{i+1}$ but $a_i\neq a_{i+1}'$. Note now that 
$$g_{i-1}a_ig_ia_{i+1}g_{i+1}\ap g_{i-1}a_ih_ia_{i+1}g_{i+1}$$
by Lemma \ref{l51a}.

We define
$$v=\left\{\begin{array}{ll}
g_0a_1\cdots g_{i-1}a_{i+2}g_{i+2}\cdots g_n &\mbox{ if } (g_{i-1},a_i)=(g_{i+1},a_{i+1}')\,,\\ [.2cm]
g_0\cdots g_{i-1}a_ih_ia_{i+1}g_{i+1}\cdots g_n\quad &\mbox{ otherwise}\,. 
\end{array}\right.$$
Note that $v$ is a landscape again and by the observations made in the previous paragraph, $u\ap v$. We say that $v$ is obtained from $u$ by \emph{uplifting the river} $g_i$, and write $u\to v$. If we need to identify the river uplifted we write instead $u\xr{g_i}v$. We use $\xr{*}$ for the reflexive and transitive closure of $\to$. Thus $\xr{*}$ is contained in $\rho$ and we have proved the following lemma:

\begin{lem}\label{to}
If $u\xr{*} v$, then $[u]=[v]$.
\end{lem}

\begin{figure}[h]
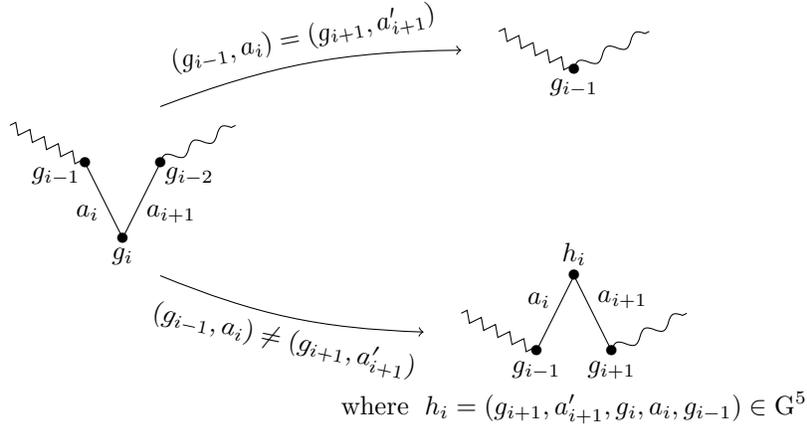

$$\tikz[scale=.5]{
\coordinate (1) at (3,3);
\coordinate (2) at (2,5);
\coordinate (3) at (0,6);
\coordinate (4) at (4,5);
\coordinate (5) at (6,6);
\draw (1) node {\small$\bullet$} node [below] {\small$g_i$};
\draw (2) node {\small$\bullet$} node [below left=-2pt] {\small$g_{i-1}$};
\draw (4) node {\small$\bullet$} node [below right=-2pt] {\small$g_{i-2}$};
\draw (2) to node[below left=-2pt] {\small$a_i$} (1);
\draw (1) to node[below right=-2pt] {\small$a_{i+1}$} (4);
\draw[snake,segment amplitude=2pt,segment length=6pt] (3)--(2);
\draw[snake=coil,segment aspect=0,segment amplitude=2pt,segment length=12pt] (4)--(5);
\coordinate (6) at (15,7.5);
\coordinate (7) at (13,8.5);
\coordinate (8) at (17,8.5);
\draw (6) node {\small$\bullet$} node [below] {\small$g_{i-1}$};
\draw[snake,segment amplitude=2pt,segment length=6pt] (7)--(6);
\draw[snake=coil,segment aspect=0,segment amplitude=2pt,segment length=12pt] (6)--(8);
\draw[->, bend left=10] (4,6.5) to node [sloped, above] {\small$(g_{i-1},a_i)=(g_{i+1},a_{i+1}')$} (12,8);
\coordinate (9) at (12,1);
\coordinate (10) at (14,0);
\coordinate (11) at (15,2);
\coordinate (12) at (16,0);
\coordinate (13) at (18,1);
\draw (10) node {\small$\bullet$} node [below] {\small$g_{i-1}$};
\draw (11) node {\small$\bullet$} node [above] {\small$h_i$};
\draw (12) node {\small$\bullet$} node [below] {\small$g_{i+1}$};
\draw (10) to node[above left=-2pt] {\small$a_i$} (11);
\draw (11) to node[above right=-2pt] {\small$a_{i+1}$} (12);
\draw[snake,segment amplitude=2pt,segment length=6pt] (9)--(10);
\draw[snake=coil,segment aspect=0,segment amplitude=2pt,segment length=12pt] (12)--(13);
\draw (15,-1.5) node {\small where $\;h_i=(g_{i+1},a_{i+1}',g_i,a_i,g_{i-1})\in \G^5$};
\draw[->, bend right=10] (4,2) to node [sloped, below] {\small$(g_{i-1},a_i)\neq(g_{i+1},a_{i+1}')$} (11,0.5);
}$$	
\caption{Illustration of the uplifting of a river.}\label{fig4}
\end{figure}

We use again a figure of \cite{LO23} to help visualize the operation of uplifting of rivers (Figure \ref{fig4}). This operation either decreases the length of a landscape or increases the height of one of its letters (keeping its length). Since the first and last letters remain the same, we conclude that after applying some uplifting of rivers to a landscape we must stop with a landscape with no rivers and with the same first and last letters. In terms of the letters subsequence, the uplifting of a river either decreases its length by two units or replaces one of its letters by another letter of height plus two. For the anchors subsequence, either it remains the same or a subword $aa'$ with $a\in A$ is deleted from it.

If $u$ is a mountain range, we always end up with a mountain if we apply uplifting of rivers while there are rivers. As explained in \cite{LO23}, we always get the same mountain independently of the order we choose to apply the upliftings of rivers. Thus, we denote by $\be_2(u)$ the only mountain such that $u\xr{*}\be(u)$. Note that $\be_2$ is only defined for mountain ranges. Let $\be(v)=\be_2(\be_1(v))$ for every $v\in\G^+$.

\begin{cor}\label{m}
For each $v\in\G^+$, $v\ap\be(v)$ and $\be(v)\in\mxm$.
\end{cor} 

\begin{proof}
This is an obvious consequence of Lemmas \ref{l52a} and \ref{to}.
\end{proof}

In the next subsection we solve the word problem for $\pres$ by showing that $\be(u)$ is the only mountain in $[u]$.

\subsection{The word problem for $\langle\G(X),\R(X)\rangle$}\label{sub46}

Let
$$\rho_1=\{(u,v)\in\G^+\times\G^+\,|\;\be(u)=\be(v)\}\,.$$
Showing that $\be(u)$ is the only mountain of $[u]$ is equivalent to prove that $\rho=\rho_1$. We already know that $\rho_1\subseteq\rho$ by Corollary \ref{m}. For the converse we need to show that both $\rho_e$ and $\rho_s$ are contained in $\rho_1$. We have already done that for $|X|=1$ in \cite[Subsection 3.6]{LO23}. The number of elements of $X$ is irrelevant for the proof of $\rho_e\subseteq\rho_1$ presented in \cite{LO23}. However, that is not the case for $\rho_s\subseteq\rho_1$, although most of the proof remains the same. Thus, we will redo only partially the proof of $\rho_s\subseteq\rho_1$, explaining what remains valid from \cite{LO23} and what needs to be added.

\begin{lem}\label{valelem}
Let $g\in \G'$ and let $u=g_0a_1g_1\cdots a_ng_n\in \lx$. Then: 
\begin{itemize}
\item[$(i)$] If $\up(g)\geq 2$ and $s\in \{l,r\}$, then 
$\la_r(g^s)*\be_1(g^{sa})*\la_l(g)\xrightarrow{*}g^sg^{sa}g\;$ and $\;\la_r(g)*\be_1((g^{sa})')*\la_l(g^s)\xr{*}g(g^{sa})'g^s\,.$
\item[$(ii)$] $\be_1(u)\xr{*} \la_l(g_0)*u*\la_r(g_n)$.
\item[$(iii)$] $\be(u)=\be_2(u)$ if $u\in \mrx$.
\item[$(iv)$] $\be(u)=u$ if $u\in \mxm$.
\end{itemize}
\end{lem}

\begin{proof}
This result restates \cite[Lemma 3.9]{LO23} for $X$ of any size. The proofs of statements $(ii)$, $(iii)$ and $(iv)$ presented in \cite{LO23} remain valid for $|X|>1$. However, that is not the case for the proof of $(i)$. In \cite{LO23}, $(i)$ is proved by induction on $\up(g)$. The number of elements of $X$ is irrelevant for the proof of the inductive step. For the base case ($\up(g)=2$), the size of $X$ matters, however, since $\G_{2,d}$ is empty for $|X|=1$, but nonempty for $|X|>1$. Hence, we need to complete the proof of $(i)$ presented in \cite{LO23} by showing that the statement of $(i)$ is valid for the elements of $\G_{2,d}$ when $|X|>1$.
	
Let $g=(g_{xx'},a,1,b,g_{yy'})\in\G_{2,d}$ for $x$ and $y$ distinct elements of $X$. Then $a\in\{1,x\}$ and $b\in\{1,y\}$. If $s=l$, then
$$\begin{array}{ll}
\la_r(g^s)*\be_1(g^{sa})*\la_l(g)& =\la_r(g_{xx'})*\be_1(a)*\la_l(g)\\ [.2cm]
&=\left\{\begin{array}{ll}
g_{xx'}111g_{xx'}1g & \mbox{ if } a=1 \\ [.2cm]
g_{xx'}111g_{xx'}x1x'g_{xx'}xg & \mbox{ if } a=x
\end{array}\right. 
\end{array}$$
In the former case, $g_{xx'}111g_{xx'}1g\xr{1} g_{xx'}1g=g^sg^{sa}g$. In the latter case, 
$g_{xx'}111g_{xx'}x1x'g_{xx'}xg\xr{*}g_{xx'}xg=g^sg^{sa}g$ by uplifting the river $1$ twice. Thus $\la_r(g^s)*\be_1(g^{sa})*\la_l(g)\xrightarrow{*}g^sg^{sa}g\;$ for $s=l$.

If $s=r$, then 
$$\begin{array}{ll}
\la_r(g^s)*\be_1(g^{sa})*\la_l(g)& =\la_r(g_{yy'})*\be_1(b)*\la_l(g)\\ [.2cm]
&=\left\{\begin{array}{ll}
g_{yy'}11a'g_{xx'}ag & \mbox{ if } b=1 \\ [.2cm]
g_{yy'}111g_{yy'}y1a'g_{xx'}ag & \mbox{ if } b=y
\end{array}\right.
\end{array}$$
In the former case, $g_{yy'}11a'g_{xx'}ag\xr{1} g_{yy'}1ga'g_{xx'}ag \xr{g_{xx'}} g_{yy'}1g=g^sg^{sa}g$. In the latter case, 
$$g_{yy'}111g_{yy'}y1a'g_{xx'}ag\xr{1,1}g_{yy'}yga'g_{xx'}ag\xr{g_{xx'}} g_{yy'}yg=g^sg^{sa}g\,.$$ 
Thus $\la_r(g^s)*\be_1(g^{sa})*\la_l(g)\xrightarrow{*}g^sg^{sa}g\;$ also for $s=r$.

The proof of $\;\la_r(g)*\be_1((g^{sa})')*\la_l(g^s)\xr{*}g(g^{sa})'g^s\;$ is similar to the proof of $\;\la_r(g^s)*\be_1(g^{sa})*\la_l(g)\xrightarrow{*}g^sg^{sa}g\;$ presented above. Hence, we can conclude that the statement of $(i)$ is valid also for the elements of $\G_{2,d}$ if $|X|>1$. We have now completed the proof of $(i)$.
\end{proof}

\begin{lem}
$\rho_s\subseteq\rho_1$.
\end{lem}

\begin{proof}
In this proof we use several times Lemma \ref{valelem}.$(ii)$. Note that 
$$\be_1(g^cg^Lg) \xr{*} \la_l(g^c)*(g^cg^Lg)*\la_r(g)=\la_l(g)*\la_r(g)=\be_1(g)\,.$$ 
Thus $\be(g^cg^Lg)=\be(g)$ and similarly $\be(gg^Rg^c)=\be(g)$. Now
$$\begin{array}{ll}
\be_1(g^r\cdot g^c\cdot g^l)\hspace*{-.1cm}& \xr{*}\;\la_l(g^r)*(g^r(g^{r_aa})'g^cg^{l_aa}g^l)*\la_r(g^l) \\ [.2cm]
& \xr{g^c}\; \la_l(g^r)*(g^r(g^{r_aa})'gg^{l_aa}g^l)*\la_r(g^l) \\ [.2cm]
& =\; \la_l(g^r)*(g^rg^{ra}g(g^{la})'g^l)*\la_r(g^l) \\ [.2cm]
& =\; \la_l(g^r)*(g^r\cdot g\cdot g^l)*\la_r(g^l)\,;
\end{array}$$
and so $\be_1(g^r\cdot g^c\cdot g^l)\ap\be_1(g^r\cdot g\cdot g^l)$. By definition of $\be_1$, we must have also $\be_1(g^Rg^cg^L)\ap\be_1(g^R g g^L)$; whence $\be(g^Rg^c g^L)=\be(g^Rgg^L)$. Note that we have shown that $\rho_s\subseteq\rho_1$.
\end{proof}

Next, we state the main result of this subsection and its proof is now obvious from the discussion yield above.

\begin{prop}\label{mrhoinv}
Let $u,v\in \Gp$. Then $\be(u)=\be(v)$ if and only if $[u]=[v]$. So, the word problem for $\pres$ is decidable.
\end{prop}

In the following section we study the structure of $\fx$. The main result is the proof that $\fx$ is a regular semigroup weakly generated by $X$.

\section{The semigroup $\fx$}\label{sub47}

This section follows the structure of \cite[Section 4]{LO23}. In fact, except for Proposition \ref{model} and the main result of this section, the proofs of the results presented in \cite[Section 4]{LO23} for $|X|=1$ are valid for any set $X$. Thus, for the sake of completeness, we restate those results here for any set $X$, but refrain ourselves from including proofs. This allows also the reader to follow easily the reasoning that leads to the main result. About Proposition \ref{model}, it corresponds to \cite[Proposition 3.12]{LO23}, a result not in Section 4 of that paper. Since the proof of \cite[Proposition 3.12]{LO23} makes use of other results whose numbering has no direct correspondence in this paper, we decided to include a proof of Proposition \ref{model} here.

The next lemma is needed before the proof that $\fox$ is regular.

\begin{lem}\cite[Lemma 4.1]{LO23}\label{updownhill}
Let $v=g_0a_1g_1\cdots a_ng_n$ be a downhill of $\G^+$. Then $v\cev{v}\ap v*\cev{v}\ap g_0$.
\end{lem}

It is now straightforward to check that $[\cev{u}]$ is an inverse of $[u]$ for every mountain $u$, thus concluding that $\fox$ is regular.

\begin{prop}\cite[Proposition 4.2]{LO23}
The monoid $\fox$ is regular and $[\cev{u}]$ is an inverse of $[u]$ for all $u\in \ltp$.
\end{prop}

We introduce a model for $\fox$ next. We will use this model to easily work with $\fox$. Let $u_1,u_2\in\mxm(X)$ and define 
$$u_1\odot u_2= \be(u_1*u_2)=\be_2(u_1*u_2)\in\mxm\,.$$
Thus $\odot$ is a binary operation on $\mxm$.

\begin{prop}\label{model}
$(\mxm(X),\odot)$ and $(\mx(X),\odot)$ are models for $\fox$ and $\fx$, respectively. Moreover, if $w=u\odot v$ for $u,v,w\in \mxm(X)$, then
\begin{itemize}
\item[$(i)$] $\la_l(u)$ is a prefix of $\la_l(w)$, while $\la_r(v)$ is a suffix of $\la_r(w)$.
\item[$(ii)$] $\up(w)\geq\max\{\up(u),\up(v)\}$, and $\up(w)=\up(u)$ {\rm [}$\up(w)=\up(v)${\rm ]} \iff\ $\ka(w)=\ka(u)$ {\rm [}$\ka(w)=\ka(v)${\rm ]}.
\end{itemize}
\end{prop}

\begin{proof}
The mapping $[u]\mapsto\be(u)$ is a well defined bijection from $\fx$ onto $\mxm(X)$ due to  Proposition \ref{mrhoinv}. But since $uv\ap\be(u)*\be(v)$, we have also
$$\be(uv)=\be(\be(u)*\be(v))=\be(u)\odot\be(v)\,,$$
that is, this mapping is an isomorphism. Note that uplifting rivers does not change the prefix $\la_l(u)$ and the suffix $\la_r(u)$ of a mountain range $u$. Then $(i)$ follows from this observation, and $(ii)$ is an obvious consequence of $(i)$.
\end{proof}

The \emph{ground} $\ep(g)$ of a letter $g\in \G'$ is defined recursively as follows: 
$$\ep(1)=\{1\}\quad\mbox{ and }\quad\ep(g)=\ep(g^l)\cup\{g\} \cup\ep(g^r)$$ 
if $\up(g)\geq 1$. Thus each $g_1\in\ep(g)\setminus\{g\}$ has height less than $g$. We write $h\preceq g$ if $h\in\ep(g)$. We define the ground of a landscape $u=g_0a_1g_1\cdots a_ng_n$ as $\ep(u)=\cup_{i=0}^n\ep(g_i)$. The next result is a technical lemma important to characterize the Green's relations on $\mxm$, especially statement $(ii)$.

\begin{lem}\cite[Lemma 4.3]{LO23}\label{ground}
Let $g,h\in \G'$ and $u,v,w\in\mxm$ such that $w=u\odot v$. Then:
\begin{itemize}
\item[$(i)$] $h\in\ep(g)$ \iff\ $\ep(h)\subseteq\ep(g)$.
\item[$(ii)$] If $h\in\ep(g)\setminus\{g\}$, then there exists $u=g_0a_1g_1\cdots a_ng_n\in \lop$ such that $n=\up(g)-\up(h)$, $g_i\in\ep(g)$ for all $0\leq i\leq n$, $g_0=h$ and $g_n=g$.
\item[$(iii)$] $\ep(u)\cup\ep(v)\subseteq\ep(w)$.
\end{itemize}
\end{lem} 

The relation $\preceq$ is clearly a partial order on $\G'$ due to statement $(i)$. Further, $\ep(u)$ is the union of the grounds of its ridges. In particular, $\ep(u)=\ep(\ka(u))$ for any mountain $u$.
 
The next proposition and the following two corollaries describe the partial orders $\leq_{\Rc}$, $\leq_{\Lc}$ and $\leq_{\Jc}$ and the Green's relations on $(M^1,\odot)$.

\begin{prop}\cite[Proposition 4.4]{LO23}\label{desR}
Let $u,v\in \mxm$. Then:
\begin{itemize}
\item[$(i)$] $u\leq_{\Rc} v$ \iff\ $\la_l(v)$ is a prefix of $\la_l(u)$. 
\item[$(ii)$] $u\leq_{\Lc} v$ \iff\ $\la_r(v)$ is a suffix of $\la_r(u)$. 
\item[$(iii)$] $u\leq_{\Jc} v$ \iff\ $\ka(v)\preceq\ka(u)$. 
\end{itemize}
\end{prop}

\begin{cor}\cite[Corollary 4.5]{LO23}
Let $u,v\in \mxm$. Then:
\begin{itemize}
\item[$(i)$]  $v$ covers $u$ for $\leq_{\Rc}$ \iff\ $\la_l(u)=\la_l(v)a\ka(u)$ for some $a\in A$.
\item[$(ii)$] $v$ covers $u$ for $\leq_{\Lc}$ \iff\ $\la_r(u)=\ka(u)a\la_r(v)$ for some $a\in A$.
\item[$(iii)$] $v$ covers $u$ for $\leq_{\Jc}$ \iff\ $\ka(v)\in\{(\ka(u))^l,(\ka(u))^r\}$.
\end{itemize}
\end{cor}

\begin{cor}\cite[Corollary 4.6]{LO23}\label{R}
Let $u,v\in\mxm$. Then:
\begin{itemize}
\item[$(i)$] $u\Rc v$ \iff\ $\la_l(u)=\la_l(v)$.
\item[$(ii)$] $u\Lc v$ \iff\ $\la_r(u)=\la_r(v)$.
\item[$(iii)$] $u\Jc v$ \iff\ $\ka(u)=\ka(v)$.
\item[$(iv)$] $u\Hc v$ \iff\ $u=v$.
\item[$(v)$] $\Dc =\Jc$.
\end{itemize}
\end{cor}

The statements $(iii)$ and $(v)$ indicate there is a natural bijection between the set of $\Dc$-classes of $\fox$ and the set $\G'$. Each $\Dc$-class of $\mxm$ is constituted by the mountains with a specific peak from $\G'$. Using the previous corollary and Lemma \ref{number} we obtain immediately the size of each $\Rc$, $\Lc$ and $\Dc$-class:

\begin{cor}\cite[Corollary 4.7]{LO23}
If $g\in \G_n$ for $n\geq 1$, then $|\Rcc_{[g]}|=2^n=|\Lcc_{[g]}|$ and $|\Dcc_{[g]}|=2^{2n}$. Further, $\Dcc_{[g]}$ has $2^n$ $\Rc$-classes and $2^n$ $\Lc$-classes. 
\end{cor}

A valley $w$ with $\si(w)=\tau(w)$ is called a \emph{canyon}. If further $w\xr{*}\si(w)$, the canyon $w$ is called a \emph{gorge}. The next proposition describes the idempotents, the inverses and the natural partial order on $\mxm$ using gorges.

\begin{prop}\cite[Proposition 4.10]{LO23}
Let $u,v\in\mxm$.
\begin{itemize}
\item[$(i)$] $u$ is an idempotent \iff\ the canyon $\la_r(u)*\la_l(u)$ is a gorge.
\item[$(ii)$] $v$ is an inverse of $u$ \iff\ $\la_r(u)*\la_l(v)$ and $\la_r(v)*\la_l(u)$ are both gorges.
\item[$(iii)$] $v<u$ \iff\ $\la_l(v)=\la_l(u)\,a_1\,u_1$ and $\la_r(v)=u_2\,a_2\, \la_r(u)$ for some $u_1\in\lop\cup\G^5$, $u_2\in\lom\cup\G^5$ and $a_1,a_2\in A$ such that $u_2\,a_2\,\ka(u)\,a_1\,u_1$ is a gorge.
\end{itemize}
\end{prop}

Note that it is hard to identify the canyons that are gorges. Thus, the previous result gives us just a theoretical description for the idempotents, the inverses and the natural partial order, but not a so useful description for practical purposes.

The next result states that the sandwich sets $S([g^Rg^c],[g^cg^L])$ are singleton sets. This is a crucial result in order to prove that $\fx$ is weakly generated by $X$.

\begin{lem}\cite[Lemma 4.8]{LO23}\label{sandwich}
For each $g\in \G^5$ with $\up(g)\geq 2$, 
$$S([g^Rg^c],[g^cg^L])=\{[g]\}$$ 
in $\fox$.
\end{lem}

Finally, we reach to the main result of this section: $\fx$ is a regular semigroup weakly generated by $\{[x]\,|\;x\in X\}$. Since $\{[g]\,|\;g\in \G\setminus\{1\}\}$ generates $\fx$, we just need to show that every regular subsemigroup $S$ of $\fx$ containing $\{[x]\,|\;x\in X\}$ contains also $\{[g]\,|\;g\in \G\setminus\{1\}\}$. In Proposition 4.9 of \cite{LO23} we proved that the previous statement is true for $|X|=1$ by induction on $\up(g)$. That proof is not complete if $X$ has more than one element. The reason is because the initial case is now more complex: we need to conclude that $\G_1$, $\G_2$ and $\{[x']\,|\;x\in X\}$ are contained in $S$, something that is trivial for $|X|=1$. The proof of the inductive step, however, continues valid. Therefore, below, we just prove that the three sets mentioned above belong to $S$, if $S$ is a subsemigroup of $\fx$ containing $\{[x]\,|\;x\in X\}$.

\begin{prop}
The regular semigroup $\fx$ is weakly generated by $\{[x]\,|\;x\in X\}$.
\end{prop}

\begin{proof}
Let $S$ be a regular subsemigroup of $\fx$ containing $\{[x]\,|\;x\in X\}$ and let $x\in X$. Since
$$\Dcc_{[x]}=\{[x],[x'],[g_{xx'}],[x'x]\}$$
by Corollary \ref{R} and $[x']$ is the only inverse of $[x]$, $S$ contains also $[x']$ and $[g_{xx'}]$. Thus $\G_1$ and $\{[x']\,|\;x\in X\}$ are contained in $S$. Let now $g\in\G_2$. Then $g^c=1$ and $S([g^Rg^c],[g^cg^L])= S([g^R],[g^L])$ for $g^R\in\{g_{xx'},[x'x]\}$ and $g^L\in\{g_{yy'},[y'y]\}$, where $x,y\in X$. Thus both $g^R$ and $g^L$ belong to $S$. Since $[g]$ is the only element of $S([g^R],[g^L])$ and $S$ is regular, $S$ must contain $[g]$; whence $\G_2\subseteq S$. As observed above, using the proof of \cite[Proposition 4.9]{LO23}, this is enough to conclude the proof of this proposition.
\end{proof}

In the next section we prove that all regular semigroups weakly generated by $X$ are homomorphic images of $\fx$ and that $\fx$ is the only (up to isomorphism) regular semigroup with this property.

\section{The universal property of $\fx$}\label{sec5}
 
A skeleton mapping is a special mapping $\phi:\G\to S^1$ where $S$ is a regular semigroup and $S^1$ is the monoid obtained from $S$ by adding an identity element $1$ if necessary. We need the following notation before we can advance and introduce the properties that define these mappings: for each $g\in\G^5$ with $\up(g)\geq 2$,
$$g^{\phi,l}=(g^c\phi)((g^{la})'\phi)(g^l\phi)(g^{la}\phi)$$ 
and
$$g^{\phi,r}=((g^{ra})'\phi)(g^r\phi)(g^{ra}\phi)(g^c\phi)\,.$$
We will see later that the skeleton mapping's properties (not yet defined) will imply that both $g^{\phi,l}$ and $g^{\phi,r}$ are idempotents of $S$. Since we will need to refer to the sandwich set $S\big(g^{\phi,r},g^{\phi,l}\big)$ in the definition of skeleton mapping, we are considering, for the moment, the non-idempotent version of the definition of these sets. Of course, after we show that $g^{\phi,l}$ and $g^{\phi,r}$ are indeed idempotents, we can then look to these sets with their usual definition for idempotents.

A \emph{skeleton mapping} (on $X$) is a mapping $\phi:\G\to S^1$, where $S$ is a regular semigroup, such that
\begin{itemize}
\item[$(i)$] $x\phi\in S$, $x'\phi\in V(x\phi)$ and $g_{xx'}\phi=(x\phi)(x'\phi)$ for any $x\in X$; 
\item[$(ii)$] $(1\phi)(a\phi)=a\phi=(a\phi)(1\phi)$ for any $a\in A$;
\item[$(iii)$] $g\phi\in S\big(g^{\phi,r},g^{\phi,l}\big)$ for any $g\in \G_i$ with $i\geq 2$.
\end{itemize}
Note that $1\phi\in E(S^1)$ since we can take $a=1$ in $(ii)$. Further $\G^5\phi\subseteq E(S)$:  $g_{xx'}\phi\in E(S)$ for any $x\in X$ by $(i)$ and $g\phi\in E(S)$ for any $g\in\G_i$ with $i\geq 2$ by $(iii)$. 

Let $s\in V(g^{\phi,l})$ and $t\in V(g^{\phi,r})$. Then $S\big(g^{\phi,r},g^{\phi,l}\big)=S\big(tg^{\phi,r}, g^{\phi,l}s\big)$ and 
$$g\phi=g^{\phi,l}s(g\phi)=(g\phi)tg^{\phi,r}\,.$$
So $(g^c\phi)(g\phi)=(g\phi)=(g\phi)(g^c\phi)$. It is now obvious by induction that $(1\phi)(g\phi)=g\phi=(g\phi)(1\phi)$ for all $g\in\G$. In these last two observations we needed to use the inverses $s$ and $t$ since we haven't prove yet that $g^{\phi,l}$ and $g^{\phi,r}$ are idempotents. After the following result there is no need to continue to use $s$ and $t$.

\begin{lem}\label{l51}
Let $S$ be a regular semigroup and $\phi:\G\to S^1$ be a skeleton mapping. Then $g^{\phi,l}$ and $g^{\phi,r}$ are idempotents of $S$ for all $g\in\G_i$ with $i\geq 2$.
\end{lem}

\begin{proof}
This result extends \cite[Lemma 5.1]{LO23} for $|X|>1$. In \cite{LO23} we proved Lemma 5.1 by induction on $i$. Once again, the proof of the inductive step remains valid for $X$ of any size. Hence, we only need to complete the base case since $\G_{2,d}$ is empty if $|X|=1$, but nonempty otherwise. Therefore, we need to verify that $g^{\phi,l}$ and $g^{\phi,r}$ are idempotents for any $g\in\G_{2,d}$ if $|X|\geq 2$.

Let $g=(g_{xx'},a,1,b,g_{yy'})\in\G_{2,d}$ with $x$ and $y$ distinct elements of $X$. Then $a\in\{1,x\}$ and $b\in\{1,y\}$. If $s=x\phi$ and $s'=x'\phi$, then $g^{\phi,l}\in\{ss',s's\}$ by $(i)$ and $(ii)$. Thus $g^{\phi,l}$ is an idempotent of $S$. Similarly, $g^{\phi,r}$ is an idempotent too.
\end{proof}

Let $\phi$ be a mapping from $X$ to a regular semigroup $S$. We can naturally extend $\phi$ to a skeleton mapping $\phi:\G\to S^1$ (we are using the same notation here for the two mappings). Indeed, for each $x\in X$, we can choose an inverse of $x\phi$ for the image of $x'$, and set $g_{xx'}\phi=(x\phi)(x'\phi)$. We can set also $1\phi=1$. Then, recursively, for each $g\in\G_i$ with $i\geq 2$, we choose an idempotent $e\in S\big(g^{\phi,r},g^{\phi,l} \big)$ and set $g\phi=e$. We obtain a skeleton mapping with this procedure. Of course, we may have several distinct skeleton mappings extending the same mapping form $X$ to $S$. The advantage of having a skeleton mapping $\phi:\G\to S^1$ is that it always has a unique extension into a homomorphism from $\fox$ to $S^1$. This is part of the next result.

\begin{prop}\cite[Proposition 5.2]{LO23}\label{r52}
If $\phi:\G\to S^1$ is a skeleton mapping, then there is a unique (semigroup) homomorphism $\varphi:\fox\to S^1$ extending $\phi$ (that is, such that $g\phi=[g]\varphi$ for all $g\in\G$). Furthermore, $\langle\G\phi\rangle=\fox\varphi$ is a regular subsemigroup of $S^1$ and a monoid with identity element $1\phi$.
\end{prop}

\begin{proof}
There is an aspect of the proof of \cite[Proposition 5.2]{LO23} that needs to be clarified first so that we can say it continues to be valid for this paper. The existence of $\varphi$ is guaranteed by showing that the natural homomorphism $\phi_1:\G^+\to S^1$ that extends $\phi$ contains $\rho_e\cup\rho_s$ in its kernel. However, the definition of $\rho_s$ here and in \cite{LO23} is not exactly the same: the pairs $(g^Rgg^L,g^Rg^cg^L)$ are replaced by $(g^r\cdot g\cdot g^l,g^r\cdot g^c\cdot g^l)$ in \cite{LO23}. Thus, from the computations yield in \cite{LO23}, we can only conclude that the pairs $(g^r\cdot g\cdot g^l,g^r\cdot g^c\cdot g^l)$ are contained in the kernel of $\phi_1$. But it is now evident that also the pairs $(g^Rgg^L,g^Rg^cg^L)$ belong to the kernel of $\phi_1$ since
$$g^Rgg^L=(g^{ra})'\,(g^r\cdot g\cdot g^l)\,g^{la}\quad\mbox{ and }\quad g^Rg^cg^L=(g^{ra})'\,(g^r\cdot g^c\cdot g^l)\,g^{la}\,.$$
Apart from this small detail, the proof presented in \cite[Proposition 5.2]{LO23} is valid here also.
\end{proof}

In the previous result, $\fox\varphi$ is a regular subsemigroup of $S^1$ and a monoid with identity $1\phi$. However, it may not be a submonoid of $S^1$: $1\phi$ is not necessarily the identity element of $S^1$.

A \emph{skeleton} of $S$ is the image of $\G\setminus\{1\}$ under some skeleton mapping. By the second part of the previous proposition, the subsemigroup of $S$ generated by a skeleton is always regular. We have justified the following corollary.

\begin{cor}
If $A$ is a skeleton of $S$, then $\langle A\rangle$ is a regular subsemigroup of $S$.
\end{cor}

We can finally prove the main result of this paper. To turn the language less cumbersome, in the following results we identify the elements of $X$ with their corresponding congruence classes in $\fx$.

\begin{teor}
All regular semigroups weakly generated by $X$ are homomorphic images of $\fx$ (under a homomorphism whose restriction to $X$ is the identity map).
\end{teor}

\begin{proof}
Let $S$ be a regular semigroup weakly generated by $X$. As explained before, we can extend the embedding of $X$ into $S$ to a skeleton mapping $\phi:\G\to S^1$, and then to a homomorphism $\varphi:\fox\to S^1$ by Proposition \ref{r52}. Furthermore $\fx\varphi$ is contained in $S$. Thus $\fx\varphi$ is a regular subsemigroup of $S$ containing $X$, and so $S=\fx\varphi$ since $S$ is weakly generated by $X$.
\end{proof}

Although the previous result tells us that any regular semigroup weakly generated by $X$ is a homomorphic image of $\fx$, we must alert that the converse is not true. There are homomorphic images of $\fx$ that are not weakly generated by $X$. For example, in \cite{LO23} we constructed a homomorphic image of $\fw(\{x\})$ that is not weakly generated by $\{x\}$. This also means that we cannot guarantee that the regular subsemigroup generated by a skeleton is indeed weakly generated by the corresponding set $X$.

We end this paper showing that any regular semigroup $S$ weakly generated by $X$ with the property that all other regular semigroups weakly generated by $X$ are homomorphic images of $S$, is isomorphic to $\fx$.

\begin{prop}
Let $S$ be a regular semigroup weakly generated by $X$ such that any other regular semigroup weakly generated by $X$ is a homomorphic image of $S$ (under a homomorphism whose restriction to $X$ is the identity map). Then $S$ is isomorphic to $\fx$ (under an isomorphism whose restriction to $X$ is the identity map).
\end{prop}

\begin{proof}
Since both $\fx$ and $S$ are weakly generated by $X$, there are surjective homomorphisms $\varphi:\fx\to S$ and $\psi:S\to\fx$ whose restrictions to $X$ are the identity map. Then $\varphi\psi:\fx\to\fx$ is also a surjective homomorphism whose restriction to $X$ is the identity map. It is enough to show now that $\varphi\phi$ is the identity map since it follows that both $\varphi$ and $\psi$ are one-to-one mappings, that is, $\varphi$ and $\psi$ are mutually inverse isomorphisms.

Clearly $[g]\varphi\psi=[g]$ for any $g\in\G_1$ since $[x']$ is the only inverse of $[x]$ in $\fx$. Since the homomorphisms between semigroups respect the sandwich sets, we conclude that $[g]\varphi\psi=[g]$ for any $[g]\in\G\setminus\{1\}$ due to Lemma \ref{sandwich}. Finally, $\varphi\psi$ is the identity isomorphism because $\fx$ is generated by $\G\setminus\{1\}$.
\end{proof}

\vspace*{.5cm}

\noindent{\bf Acknowledgments}: This work was partially supported by CMUP, member of LASI, which is financed by (Portuguese) national funds through FCT - Fundação para a Ciência e a Tecnologia, I.P., under the project with reference UIDB/00144/2020. The author also would like to acknowledge the importance of the GAP software \cite{gap}, and its Semigroup package \cite{mitchell}, in the research presented in this paper: the several simulations computed in GAP allowed the emergence of the pattern used to construct the set $\G$.

\end{document}